\documentclass[12pt]{article}
\usepackage{epsfig,latexsym,amsfonts,amsmath,amssymb}

\def\textsf{\bf}  %%%%% here!!!!

\newtheorem{Theorem}{Theorem}
\newtheorem{Lemma}%[Theorem]
{Lemma}
{Proposition}
\newtheorem{Corollary}%[Theorem]
{Corollary}
\newtheorem{Remark}%[Theorem]
{Remark}

\def\proof{\noindent{\textbf{Proof. }}}

\def\SoOm{{H^s(\Omega)}}
\def\SoOmNull{H^s_0(\Omega)}
\def\tildeSobolev{\widetilde{H}^s(\Omega)}
\def\SoRn{H^s(\mathbb R^n)}

\def\irn{\int\limits_{\mathbb R^n}}
\def\div{{\rm div}}
\def\sstar{{2^*_s}}

\begin{document}

\title %{\vspace{-10mm}
{On fractional Laplacians}

\author{Roberta Musina\footnote{Dipartimento di Matematica ed Informatica, Universit\`a di Udine,
via delle Scienze, 206 -- 33100 Udine, Italy. Email: {roberta.musina@uniud.it}. 
{Partially supported by Miur-PRIN 2009WRJ3W7-001 ``Fenomeni di concentrazione e {pro\-ble\-mi} di analisi geometrica''.}}~ and
Alexander I. Nazarov\footnote{St.Petersburg Department of Steklov Institute, Fontanka, 27, St.Petersburg, 191023, Russia
and St.Petersburg State University, 
Universitetskii pr. 28, St.Petersburg, 198504, Russia. E-mail: {al.il.nazarov@gmail.com}. Supported by RFBR grant 11-01-00825 and by
St.Petersburg University grant 6.38.670.2013.}
}

\date{}

\maketitle

\footnotesize

\noindent
{\bf Abstract.} We compare two natural types of fractional Laplacians $(-\Delta)^s$, namely,
the ``Navier'' and the ``Dirichlet'' ones. We show that for $0<s<1$ their difference is
positive definite and positivity preserving. Then we prove the coincidence of the 
Sobolev constants for these two fractional Laplacians. 

\normalsize

\bigskip\bigskip

%\section{Introduction}
In recent years a lot of efforts have been invested in studying variational problems
involving nonlocal differential operators. Contrary to  the standard 
Laplacian, that acts by pointwise differentiation, these operators
 are usually defined via 
global integration and permits to describe, for instance, diffusion processes in 
presence of long range interactions.  In this context, a 
model  operator is the fractional Laplacian $(-\Delta)^s$, for
$0<s<1$.

In modeling diffusion processes for a material that is confined in a bounded
region  $\Omega$ one clearly has to take into account the nonlocal 
nature of the problem. As a matter of fact, the boundary conditions 
that naturally can be coupled to equations of the form
\begin{equation*}
\label{eq:equation}
(-\Delta)^s u=f\quad\text{in $\Omega$}
\end{equation*}
do reflect long-range interactions. Usually two types of such
boundary conditions are considered. Both of them arise {\em together}
with the fractional Laplacian operator and we call them 
{\em Navier}-type and {\em Dirichlet}-type, respectively. 

Let us first remind some well-known facts concerning  {\em poliharmonic}
operators of  order $2k$ (here $k\ge 1$ is any integer number) in a sufficiently smooth bounded domain $\Omega$. %It is well known that 
The Navier boundary conditions for the operator $(-\Delta)^k$,
are defined as follows:
$$ u\big|_{\partial\Omega}=\Delta u\big|_{\partial\Omega}=\Delta^2 u\big|_{\partial\Omega}=\dots=\Delta^{k-1} u\big|_{\partial\Omega}=0.
$$
Obviously, the corresponding operator $(-\Delta_{\Omega})^k_N$ is the $k$th power of conventional Dirichlet Laplacian in the sense of spectral theory, 
and it can be defined by its quadratic form 
$$((-\Delta_{\Omega})^k_Nu,u)=\sum\nolimits_j\lambda_j^k\cdot|(u,\varphi_j)|^2.
$$
Here, $\lambda_j$ and $\varphi_j$ are eigenvalues and eigenfunctions of the Dirichlet Laplacian in $\Omega$, respectively.

On the other hand, the Dirichlet boundary conditions for the operator $(-\Delta)^k$ are defined as follows:
$$ u\Big|_{\partial\Omega}=\frac{\partial u}{\partial {\bf n}}\Big|_{\partial\Omega}=\frac{\partial^2 u}{\partial {\bf n}^2}\Big|_{\partial\Omega}=\dots
=\frac{\partial^{k-1} u}{\partial {\bf n}^{k-1}}\Big|_{\partial\Omega}=0,
$$
where ${\bf n}$ is the unit exterior normal vector to $\partial\Omega$.
It is easy to see that the quadratic form of corresponding operator $(-\Delta_{\Omega})^k_D$ can be obtained as the restriction of the quadratic form for
the operator $(-\Delta)^k$ in $\mathbb{R}^n$ to the set of functions supported in $\Omega$:
$$((-\Delta_{\Omega})^k_Du,u)=\int\limits_{\mathbb R^n}|\xi|^{2k}|{\cal F}u(\xi)|^2 d\xi,
$$
where $\mathcal F$ is the 
Fourier transform
$$
\mathcal F{u}(\xi)=\frac{1}{(2\pi)^{n/2}}\int\limits_{\mathbb R^n} e^{-i~\!\!\xi\cdot x}u(x)~\!dx.
$$

Now for arbitrary $s>0$ we can define the ``Navier'' fractional Laplacian by the quadratic form 
$$Q_s^N[u]=((-\Delta_{\Omega})^s_Nu,u):=\sum\nolimits_j\lambda_j^s\cdot|(u,\varphi_j)|^2
$$
%with domain
%$${\rm Dom}(Q_s^N)=\{u\in L_2(\Omega)\,:\,Q_s^N[u]<\infty\}
%$$
and the ``Dirichlet'' fractional Laplacian by the quadratic form
$$
Q_s^D[u]=((-\Delta_{\Omega})^s_Du,u):=\int\limits_{\mathbb R^n}|\xi|^{2s}|{\cal F}u(\xi)|^2 d\xi
$$
with domains, respectively,
$$\aligned
{\rm Dom}(Q_s^N)=&\ \{u\in L_2(\Omega)\,:\,Q_s^N[u]<\infty\};\\
{\rm Dom}(Q_s^D)=&\ \{u\in L_2(\mathbb{R}^n)\,:\,{\rm supp}\, u\subset\overline{\Omega},\ Q_s^D[u]<\infty\}.
\endaligned
$$
For $s=1$, these two operators evidently coincide. We emphasize that, in contrast to $(-\Delta_{\Omega})^s_N$, the operator $(-\Delta_{\Omega})^s_D$
is not the $s$th power of the Dirichlet Laplacian for $s\ne1$. 

Recall that the Sobolev space $\SoRn=W^s_2(\mathbb{R}^n)$ is defined by the completion of ${\cal C}^\infty_0(\mathbb{R}^n)$ with
respect to the Hilbertian norm
$$
\|u\|_s^2=\int\limits_{\mathbb R^n}\left(1+|\xi|^2\right)^s|\mathcal Fu(\xi)|^2~\!d\xi,
$$
see for instance Section~2.3.3 of the classical monograph \cite{Tr}.
For a bounded domain $\Omega$ with Lipschitz 
boundary we put
$$
\SoOm=\left\{u\big|_{\Omega}\,:\,u\in\SoRn\right\},
$$
see \cite[Sec.~4.2.1]{Tr}  and the extension theorem in \cite[Sec.~4.2.3]{Tr}.

Also we introduce the space
$$\tildeSobolev=\{u\in \SoRn\,:\,{\rm supp}\, u\subset\overline{\Omega}\}.
$$
By Theorem 4.3.2/1 \cite{Tr}, for $s-\frac{1}{2}\notin\mathbb{Z}$ this space coincides with $\SoOmNull$ that is the closure of
${\cal C}^{\infty}_0(\Omega)$ in $\SoOm$ while for $s-\frac{1}{2}\in\mathbb{Z}$ one has 
$\tildeSobolev\subsetneq\SoOmNull$. Moreover, ${\cal C}^{\infty}_0(\Omega)$ is dense in $u\in \tildeSobolev$.

In what follows, we assume $0<s<1$. In this case both the operators $(-\Delta_{\Omega})^s_N$ and $(-\Delta_{\Omega})^s_D$
were considered in many articles on semilinear equations, see for instance
\cite{BCDS,FW,ROS1,ROS2,SV,SV2,Tan}, and compared in \cite{F}, \cite{SV3}.
We establish further relations between them.\medskip

%Our main results are stated in Theorems \ref{T:main1} and \ref{T:main2} below.
We start with a preliminary result.

%\section{The main result}
\begin{Lemma}\label{domain}
The domains of quadratic forms ${\rm Dom}(Q_s^N)$ and ${\rm Dom}(Q_s^D)$ coincide with $\tildeSobolev$.
\end{Lemma}

\proof For $Q_s^D$ the conclusion follows directly from definition. For $Q_s^N$, using the notation of interpolation spaces from \cite{Tr}, we write the following chain of equalities:
$$\aligned
{\rm Dom}(Q_s^N)&={\rm Dom}((-\Delta_{\Omega})^{s/2}_N) &&\mbox{(\cite[Theorem 10.1.1]{BS})}\\
&=\big[L_2(\Omega),{\rm Dom}((-\Delta_{\Omega})^{1/2}_N)\big]_s &&\mbox{(\cite[Theorem 1.15.3]{Tr})}\\
&=\big[L_2(\Omega),{\rm Dom}(Q_1^N)\big]_s &&\mbox{(\cite[Theorem 10.1.1]{BS})}\\
&=\big[L_2(\Omega),\widetilde H^1(\Omega)\big]_s=\tildeSobolev, &&\mbox{(\cite[Theorem 4.3.2/2]{Tr})}
\endaligned
$$
and the Lemma follows.\hfill$\square$\medskip

We point out three elementary corollaries of this lemma.

\begin{Corollary} The following relations hold in $\tildeSobolev$:
$$Q_s^D[u]\asymp \|u\|_s^2;\qquad Q_s^N[u]\asymp \|u\|_s^2.
$$
\end{Corollary}

\proof
This statement immediately follows from Lemma \ref{domain}, the Friedrichs inequality and the closed graph theorem.\hfill$\square$%\medskip

\begin{Corollary} The operators $[(-\Delta_{\Omega})^s_N]^{-1}(-\Delta_{\Omega})^s_D$ and $[(-\Delta_{\Omega})^s_D]^{-1}(-\Delta_{\Omega})^s_N$ are bounded in $\tildeSobolev$.
\end{Corollary}

\proof
Indeed, both operators $(-\Delta_{\Omega})^s_N$ and $(-\Delta_{\Omega})^s_D$ are positive definite. Thus, by Lemma 1 they map isomorphically the space $\tildeSobolev$
onto the dual space, and the statement follows.\hfill$\square$%\medskip

\begin{Corollary} The operators $(-\Delta_{\Omega})^s_N$ and $(-\Delta_{\Omega})^s_D$ have discrete spectra.
\end{Corollary}

\proof
Since $\tildeSobolev$ is compactly embedded into $L_2(\Omega)$, see, e.g., \cite[Theorem 4.10.1]{Tr}, the statement follows from \cite[Theorem 10.1.5]{BS}.\hfill$\square$\medskip

Now we prove the main results of our paper. Namely, we show that for $0<s<1$ the operator $(-\Delta_{\Omega})^s_N-(-\Delta_{\Omega})^s_D$ is 
positivity preserving (Theorem \ref{T:main1}) and positive definite (Theorem \ref{T:main2}).

\begin{Theorem} 
\label{T:main1}
%The operators $(-\Delta_{\Omega})^s_N$ and $(-\Delta_{\Omega})^s_D$ are not the same. More precisely, 
For $u\in \tildeSobolev$, $u\ge0$, the following relation holds in the sense of distributions:
\begin{equation}
(-\Delta_{\Omega})^s_Nu \ge (-\Delta_{\Omega})^s_Du.
\label{pos_pres}
\end{equation}
If $u\not\equiv0$ then (\ref{pos_pres}) holds with strict sign.
\end{Theorem}

\proof
 In the paper \cite{CaSi}, see also \cite{CT}, the fractional Laplacian in $\mathbb R^n$ was connected with the so-called 
{\it harmonic extension in $n+2-2s$ dimensions}. Namely, it was shown that the function $w_s^D(x,y)$ minimizing the weighted Dirichlet integral
$$
{\cal E}_s^D(w)=\int\limits_0^\infty\!\int\limits_{\mathbb{R}^n} y^{1-2s}|\nabla w(x,y)|^2\,dxdy
$$
over the set
$${\cal W}^D(u)=\Big\{w(x,y)\,:\,%\displaystyle\int\limits_0^\infty\!\int\limits_{\mathbb{R}^n} 
%y^{1-2s}|\nabla w(x,y)|^2\,dxdy
{\cal E}_s^D(w)<\infty~,\ \ w\big|_{y=0}=u\Big\},
$$
satisfies
%Moreover, the relation (\ref{extension_D}) immediately implies
\begin{equation}
Q_s^D[u]=\frac {C_s}{2s}\cdot {\cal E}_s^D(w_s^D),
\label{quad_D}
\end{equation}
where the constant $C_s$ is given by
$$
C_s:=\frac{4^s\Gamma(1+s)}{\Gamma(1-s)}.
$$
Moreover, $w_s^D(x,y)$ is the solution of the BVP
$$-\div (y^{1-2s}\nabla w)=0\quad \mbox{in}\quad \mathbb R^n\times\mathbb R_+;\qquad w\big|_{y=0}=u,
$$
and for sufficiently smooth $u$
\begin{equation}
(-\Delta_{\Omega})^s_Du(x)=-C_s\cdot\lim\limits_{y\to0^+} 
\frac{w_s^D(x,y)-u(x)}{y^{2s}},\qquad x\in\Omega.
\label{extension_D}
\end{equation}

In \cite{ST} this approach was generalized to quite general situation. In particular, it was shown that
the function $w_s^N(x,y)$ minimizing the energy integral
\begin{equation*}
{\cal E}_s^N(w)=\int\limits_0^\infty\!\int\limits_{\Omega} y^{1-2s}|\nabla w(x,y)|^2\,dxdy
\label{energy_N}
\end{equation*}
over the set 
$$
%{\cal W}^N_{\Omega}(u)=\{w(x,y)\,:\,w\big|_{y=0}=u,\quad w\big|_{x\in\partial\Omega}=0\},
{\cal W}^N_{\Omega}(u)=\{w(x,y)\in{\cal W}^D(u)\,:\,w\big|_{x\in\partial\Omega}=0\},
$$
satisfies
%and (\ref{extension_N}) implies
\begin{equation}
Q_s^N[u]=\frac {C_s}{2s}\cdot {\cal E}_s^N(w_s^N).
\label{quad_N}
\end{equation}
Moreover, $w_s^N(x,y)$ is the solution of the BVP
$$-\div (y^{1-2s}\nabla w)=0\quad \mbox{in}\quad \Omega\times\mathbb R_+;\qquad w\big|_{y=0}=u;
\qquad w\big|_{x\in\partial\Omega}=0,
$$
and for sufficiently smooth $u$
\begin{equation}
(-\Delta_{\Omega})^s_Nu(x)=-C_s\cdot\lim\limits_{y\to0^+} 
\frac{w_s^N(x,y)-u(x)}{y^{2s}}.
\label{extension_N}
\end{equation}

Note that formulae (\ref{extension_D}) and (\ref{extension_N}) imply
\begin{equation}
(-\Delta_{\Omega})^s_Nu-(-\Delta_{\Omega})^s_Du=C_s\cdot\lim\limits_{y\to0^+} 
\frac{w_s^D(\cdot,y)-w_s^N(\cdot,y)}{y^{2s}}
\label{differ}
\end{equation}
at least for $u\in{\cal C}^{\infty}_0(\Omega)$. However, by approximation argument the relation (\ref{differ}) holds for $u\in\tildeSobolev$ in the sense 
of distributions.

By the maximum principle the assumption $u\ge0$ implies $w_s^D\ge0$ in $\mathbb R^n\times\mathbb R_+$. Thus, $w_s^D\ge w_s^N$ at 
$\partial\Omega\times\mathbb R_+$ and, again by the maximum principle, $w_s^D\ge w_s^N$ in $\Omega\times\mathbb R_+$. Hence (\ref{differ})
gives (\ref{pos_pres}).

Let, in addition, $u\not\equiv0$. Then the strong maximum principle gives 
$$W:=w_s^D-w_s^N>0 \qquad\mbox{in}\quad \Omega\times\mathbb R_+.
$$
After changing of the variable $t=y^{2s}$ the function $W(x,t)$ solves %meets the following relations:
\begin{equation}
\Delta_xW+4s^2t^{\frac{2s-1}s}W_{tt}=0\quad\mbox{in}\quad \Omega\times\mathbb R_+;\qquad W\big|_{t=0}=0.
\label{Maz}
\end{equation}
The differential operator in (\ref{Maz}) satisfies the assumptions of \cite[Theorem 1.4]{ABMMZ} (the boundary point lemma)
at any point $(x_0,0)\in\Omega\times\{0\}$, with
$$
A(x,t)=\left(\begin{array}{cc}
	I_n&0  \\
	0&4s^2t^{\frac{2s-1}{s}}
	\end{array}\right)~,\quad
\omega({r})=r~,
\quad h=\left(\begin{array}{c} 0\\1 \end{array}\right)
~\!.
$$
In particular, the key requirement
$$\limsup\limits_{t\to 0^+}\frac{n+4s^2t^{\frac{2s-1}s}}{\frac{\widetilde\omega(t)}{t}\cdot 4s^2t^{\frac{2s-1}s}}<\infty~\!,
$$
compare with (1.18) in \cite{ABMMZ}, is satisfied by choosing
$$
\widetilde\omega(t)=\begin{cases}
t,& 0<s<\frac{1}{2\vphantom{\frac{1}{B^B}}}; \\
 t^\frac{1-s}s, & \frac{1\vphantom{\frac{B^B}{1}}}{2}\le s<1. 
 \end{cases}
$$
%{\,\min\{\frac{1-s}s,1\}}$, 
By Theorem 1.4 in \cite{ABMMZ} we have for any $x\in\Omega$
$$\liminf\limits_{y\to 0^+}\frac{W(x,y)}{y^{2s}}=\liminf\limits_{t\to 0^+}\frac{W(x,t)}{t}>0.
$$
This completes the proof in view of (\ref{differ}).\hfill$\square$%\medskip

\begin{Remark}
 In \cite{F} this fact was proved for $s= 1/2$ and for smooth $u$.
\end{Remark}

\begin{Theorem} 
\label{T:main2}
For $u\in \tildeSobolev$, the following relation holds:
\begin{equation}
\label{pos_def}
((-\Delta_{\Omega})^s_Nu,u) \ge ((-\Delta_{\Omega})^s_Du,u).
\end{equation}
If $u\not\equiv0$ then (\ref{pos_def}) holds with strict sign.
\end{Theorem}

\proof
Note that if we assume a function $w\in{\cal W}^N_{\Omega}(u)$ to be extended by zero to $(\mathbb{R}^n\setminus\Omega)\times\mathbb{R}_+$ then
evidently ${\cal W}^N_{\Omega}(u)\subset{\cal W}^D(u)$ and ${\cal E}_s^N={\cal E}_s^D\big|_{{\cal W}^N_{\Omega}(u)}$. Therefore, (\ref{quad_D}) and (\ref{quad_N}) provide
$$Q_s^N[u]=\frac {C_s}{2s}\cdot \inf\limits_{w\in{\cal W}^N_{\Omega}(u)}{\cal E}_s^N(w)
\ge \frac {C_s}{2s}\cdot \inf\limits_{w\in{\cal W}^D(u)}{\cal E}_s^D(w)=Q_s^D[u],
$$
and (\ref{pos_def}) follows.

To complete the proof, we observe that %formulae (\ref{quad_D}) and (\ref{quad_N}) hold for all $u\in\tildeSobolev$. Moreover, 
for $u\not\equiv0$ the function $w_s^N$ (extended by zero) cannot be a solution
of the homogeneous equation in the whole half-space $\mathbb R^n\times\mathbb R_+$ since such a solution 
should be analytic in the half-space. Thus, it cannot provide $\inf\limits_{w\in{\cal W}^D(u)}{\cal E}_s^D(w)$, and the last statement follows.
\hfill$\square$%\medskip

\begin{Corollary} Let us denote the eigenvalues of $(-\Delta_{\Omega})^s_N$ and $(-\Delta_{\Omega})^s_D$ by $\lambda_{s,j}^N$ and $\lambda_{s,j}^D$, respectively,
and enumerate them in ascending order according to their multiplicities. Then
$$\lambda_{s,j}^N>\lambda_{s,j}^D,\qquad j\in\mathbb{N}.
$$
\end{Corollary}

\proof
This assertion immediately follows from Theorem \ref{T:main2} by the Courant variational principle (see, e.g., \cite[Theorem 10.2.2]{BS}).\hfill$\square$%\medskip

\begin{Remark}
In \cite{ChSo} this result was obtained with $\ge$ sign by using considerably more complicated techniques. Theorem \ref{T:main2}
(also with $\ge$ sign) can be extracted from \cite[Lemma 19]{FG}.
\end{Remark}

Now we continue to compare the quadratic forms $Q_s^N$ and $Q_s^D$. First, we observe that for $u\in \tildeSobolev$ the form $Q_s^D[u]$
evidently does not change if we consider arbitrary domain $\Omega'\supset\Omega$ instead of $\Omega$. In contrast, the form $Q_s^N[u]$
depends on $\Omega'\supset\Omega$. To emphasize this dependence we introduce the notation $Q^N_s[u;\Omega]$.

\begin{Theorem}
\label{T:new}
If $u\in\tildeSobolev$, then
$$
Q^D_s[u]=\inf_{\Omega'\supset\Omega} Q^N_s[u;\Omega'],
$$
where the infimum is taken over the set of smooth bounded domains in $\mathbb R^n$.
\end{Theorem}

\proof
For $\Omega'\supset\Omega$ we have ${\cal W}^N_{\Omega}(u)\subset {\cal W}^N_{\Omega'}(u)$ for any
$u\in \tildeSobolev\subset\widetilde H^s(\Omega')$ (we recall that
$w\in{\cal W}^N_{\Omega}(u)$ are assumed extended by zero to $(\mathbb{R}^n\setminus\Omega)\times\mathbb{R}_+$).
By (\ref{quad_N}), $Q^N_s[u;\Omega]$ is monotone decreasing with respect to the domain inclusion.
Taking Theorem \ref{T:main2} into account, we obtain 
\begin{equation}
\label{eq:monotone}
Q^D_s[u]\le Q^N_s[u;\Omega']\le Q^N_s[u;\Omega].
\end{equation}

Denote by $w=w^D_s$ the Caffarelli--Silvestre extension of $u$.  Formula
(\ref{quad_D}) 
implies that the quantity
$$
\int\limits_0^\infty \frac{1}{r}\Big\{r\int\limits_0^\infty\int\limits_{{\mathbb S}_{r}} y^{1-2s}|\partial_y w(x,y)|^2
~\!d{\mathbb S}_r(x)dy\Big\} 
dr=\int\limits_0^\infty\irn y^{1-2s}|\partial_y w(x,y)|^2~\!dxdy
%<\infty,
$$
is finite (here $\mathbb S_r$ is the sphere of radius $r$ in $\mathbb R^n$).
Since the function $r\mapsto r^{-1}$
is not integrable at $\infty$, there exists a sequence $r_h\to\infty$ such that
the balls $B_{r_h}$ contain $\Omega$ and
\begin{equation}
\label{eq:rh2}
r_h\int\limits_0^\infty\int\limits_{{\mathbb S}_{r_h}} y^{1-2s}|\partial_y w(x,y)|^2
~\!d{\mathbb S}_{r_h}(x)dy \to 0.
\end{equation}

Next, for any $y\ge 0$ let $\phi_h(\cdot,y)$ be the harmonic extension of
$w(\cdot,y)$ on $B_{r_h}$, that is,
$$
%\begin{cases}
-\Delta \phi_h(\cdot,y)=0\quad\text{in} \ B_{r_h};\qquad%\\
\phi_h(\cdot,y)=w(\cdot,y)\quad\text{on} \ {\mathbb S}_{r_h}~\!.
%\end{cases}
$$
Clearly, $\phi_h(\cdot,0)\equiv 0$. 

Finally, for $x\in B_{r_h}$ and $y\ge0$ we define
$$
w_h(x,y)=
w(x,y)-\phi_h(x,y)~\!.
$$

Let us estimate ${\cal E}^N_s(w_h)$. 
We start with term involving derivatives with respect to $x$. Since $\phi(\cdot,y)$
is harmonic in $B_{r_h}$, for any fixed $y>0$ we have
$$
0=\int\limits_{B_{r_h}} (-\Delta\phi_h)(\phi_h-w)dx=\int\limits_{B_{r_h}} |\nabla\phi_h|^2dx-
\int\limits_{B_{r_h}} \nabla\phi_h\cdot\nabla w~\! dx.
$$
Therefore
\begin{multline*}
\int\limits_0^\infty\int\limits_{B_{r_h}}  y^{1-2s}|\nabla_xw_h|^2~\!dxdy\\
=\int\limits_0^\infty\int\limits_{B_{r_h}}  y^{1-2s}|\nabla_xw|^2~\!dxdy
-\int\limits_0^\infty\int\limits_{B_{r_h}}  y^{1-2s}|\nabla_x\phi_h|^2~\!dxdy\\
\le\int\limits_0^\infty\int\limits_{\mathbb{R}^n}  y^{1-2s}|\nabla_xw|^2~\!dxdy.
\end{multline*}
%As for derivative with respect to $y$ {\bf add references}, we notice that for any fixed

Next, by differentiating the Poisson formula, we notice that the function $\partial_y\phi_h(\cdot,y)$ solves
$$
-\Delta \partial_y\phi_h(\cdot,y)=0\quad\text{in} \ B_{r_h};\qquad%\\
\partial_y\phi_h(\cdot,y)=\partial_yw(\cdot,y)\quad\text{on} \ {\mathbb S}_{r_h}~\!.
$$
Therefore, $|\partial_y\phi_h(\cdot,y)|^2$ is subharmonic in $B_{r_h}$ and thus
the function
$$
\rho\mapsto\frac{1}{\rho^{n-1}}\int\limits_{{\mathbb S}_\rho}|\partial_y\phi_h(x,y)|^2~\! d{\mathbb S}_\rho(x)
$$
is nondecreasing for $\rho\in(0,r_h)$. This implies
\begin{eqnarray*}
\int\limits_{B_{r_h}}|\partial_y\phi_h(x,y)|^2~\!dx&=&
\int\limits_0^{r_h}\int\limits_{{\mathbb S}_{\rho}} |\partial_y\phi_h(x,y)|^2~\!d{\mathbb S}_\rho(x)d\rho\\
%&=&
%\int\limits_0^{r_h}\rho^{n-1}\Big\{\frac{1}{\rho^{n-1}}\int\limits_{{\mathbb S}_{\rho}} |\partial_y\phi_h(x,y)|^2
%~\!d{\mathbb S}_{\rho}(x)\Big\}d\rho\\
&\le&
\frac{1}{r_h^{n-1}}\int\limits_{{\mathbb S}_{r_h}} |\partial_y\phi_h(x,y)|^2~\!d{\mathbb S}_{r_h}(x)\cdot
\Big(\int\limits_0^{r_h}\rho^{n-1}d\rho\Big)\\
&=&
\frac{r_h}{n}\cdot\int\limits_{{\mathbb S}_{r_h}} |\partial_yw(x,y)|^2~\!d{\mathbb S}_{r_h}(x).
\end{eqnarray*}
Therefore, by (\ref{eq:rh2})
$$
\int\limits_0^\infty\int\limits_{B_{r_h}} y^{1-2s}|\partial_y\phi_h|^2~\!dxdy=o(1),\qquad r_h\to\infty,
$$
and we arrive at
\begin{eqnarray*}
\int\limits_0^\infty\int\limits_{\Omega} y^{1-2s}|\partial_yw_h|^2~\!dxdy&=&
\int\limits_0^\infty\int\limits_{B_{r_h}} y^{1-2s}|\partial_yw|^2~\!dxdy +
o(1)\\
&\le&
\int\limits_0^\infty\irn y^{1-2s}|\partial_y w|^2~\!dxdy+o(1).
\end{eqnarray*}
The above calculations imply
$$
{\cal E}^N_s(w_h)\le {\cal E}^D_s(w)+o(1),\qquad r_h\to\infty.
$$
Since $w_h(\cdot,0)=u$, we have $w_h\in {\cal W}^N_{B_{r_h}}(u)$.
Therefore, by (\ref{quad_D}) and (\ref{quad_N}) we obtain
\begin{equation}
\label{eq:tesi}
Q^N_s[u;B_{r_h}]\le \frac {C_s}{2s}\cdot {\cal E}^N_s(w_h)
\le \frac {C_s}{2s}\cdot{\cal E}^D_s(w)+o(1)= Q^D_s[u]+o(1)~\!.
\end{equation}
The conclusion readily follows by comparing (\ref{eq:monotone}) and (\ref{eq:tesi}).
\hfill$\square$%\medskip

\begin{Remark}
\label{R:new}
Assume that $0\in\Omega$ and put $\alpha\Omega=\{\alpha x\,:\,x\in\Omega\}$.
Thanks to (\ref{eq:monotone}), the proof above shows indeed that
$$
Q^D_s[u]=\lim_{\alpha\to\infty}Q^N_s[u;\alpha\Omega]
\quad\text{for any}\quad u\in\tildeSobolev.
$$
\end{Remark}

\begin{Corollary}
\label{T:almost}
$\displaystyle{\inf_{\scriptstyle v\in \tildeSobolev\atop\scriptstyle  v\ne 0}
\frac{Q^N_s[v]}
{Q^D_s[v]}=1}$.
\end{Corollary}

\proof
By Theorem \ref{T:main2} the infimum in the statement cannot be smaller than~$1$.
To conclude the proof we can suppose $0\in\Omega$.
Given $u\in \tildeSobolev$, put $u_\alpha(x)=u(\alpha x)$. Then, by homogeneity, 
$$\frac{Q^N_s[u_\alpha]}{Q^D_s[u_\alpha]}=\frac{Q^N_s[u;{\alpha\Omega}]}{Q^D_s[u]},
$$
and the statement follows  by Remark \ref{R:new}.\hfill$\square$\medskip

Finally, we assume\footnote{This is a restriction only for $n=1$.} $n> 2s$ and deal with the Sobolev embedding
$$
\tildeSobolev\hookrightarrow L_\sstar(\Omega)~,\quad \sstar=\frac{2n}{n-2s}.
$$
We introduce the Sobolev constants for the  Navier and  Dirichlet fractional Laplacians:
$$
S^N_s(\Omega)=\inf_{\scriptstyle u\in \tildeSobolev\atop\scriptstyle  u\ne 0}
\frac{Q^N_s[u]}
{\|u\|^2_{L_\sstar(\Omega)}}~,\qquad
S^D_s(\Omega)=\inf_{\scriptstyle u\in \tildeSobolev\atop\scriptstyle  u\ne 0}
\frac{Q^D_s[u]}
{\|u\|^2_{L_\sstar(\Omega)}}~\!.
$$
The constant $S^D_s(\Omega)$ is obviously invariant with respect to dilations in $\mathbb R^n$ and 
monotone increasing with respect to the domain inclusion. Thus, $S^D_s$ does not depend on the domain $\Omega$ 
(see \cite{SV4} for a more general problem). It follows from the proof of Theorem \ref{T:new} that
$S^N_s(\Omega)$ has the same property.

Next, it is evident that $S^D_s(\Omega)=S^D_s(\mathbb{R}^n)$. In \cite{CoTa} it has been shown that
the best constant 
$$
S^D_s(\mathbb{R}^n)={(4\pi)^s}~\frac{\Gamma\left(\frac{n+2s}{2}\right)}
{\Gamma\left(\frac{n-2s}{2}\right)}~\!
\left[\frac{\Gamma({n/2})}{\Gamma({n})}\right]^{2s/n}
$$
 is attained, up to dilations, translations and multiplications, only by the function
\begin{equation}\label{U}
U(x)=\left(1+|x|^2\right)^{\frac{2s-n}{2}}~\!.
\end{equation}

The equality $S^N_s(\Omega)=S^D_s(\mathbb{R}^n)$ is not so trivial. Actually, in some papers on semilinear equations
the equality between two Sobolev constants has been used but, as far as we know,
never rigorously proved. Here we fill this gap. 

\begin{Theorem}
\label{T:main3} The constant $S^N_s(\Omega)$ is equal to $S^D_s(\mathbb{R}^n)$ and is not attained on $\tildeSobolev$.
 \end{Theorem}
 
 \proof
 We assume again that $0\in\Omega$ and put $\mathcal S:=S^D_s(\mathbb R^n)$.
Thanks to (\ref{pos_def}), it holds that 
$S^N_s(\Omega)\ge S^D_s(\Omega)=\mathcal S$. So, we have to show that $S^N_s(\Omega)\le \mathcal S$. 
 
Fix any $\varepsilon>0$ and take a  nontrivial function $u\in {\cal C}^\infty_0(\mathbb R^n)$ such that
$$
{\|u\|^2_{L_\sstar(\mathbb R^n)}}(\mathcal S+\varepsilon)\ge {Q^D_s[u]}.
$$
For $\alpha>1$ the function
$u_\alpha(x):=u(\alpha x)$
belongs to $\tildeSobolev$. Therefore, by homogeneity and
thanks to Remark \ref{R:new}, we have
$$
S^N_s(\Omega) \le \frac{Q^N_s[u_\alpha;\Omega]}{\|u_\alpha\|^2_{L_\sstar(\Omega)}}\\
=\frac{Q^N_s[u;\alpha\Omega]}{\|u\|^2_{L_\sstar(\alpha\Omega)}}\\
=\frac{Q^D_s[u]+o(1)}{\|u\|^2_{L_\sstar(\mathbb R^n)}}\le
\mathcal S+\varepsilon+o(1)
$$
 as $\alpha\to\infty$. The first statement of Theorem immediately follows. The second one is obvious.\hfill$\square$

\begin{Remark}
 The statements of Theorems \ref{T:main1}--\ref{T:main3} hold also for unbounded Lipschitz domains if we define
$(-\Delta_{\Omega})^s_N$ as the $s$th power of $-\Delta_{\Omega}$ (in this case it can have continuous spectrum).
The proof runs with minor changes.
\end{Remark}

\medskip
\noindent
{\bf Acknowledgments}.
We are greatful to Rupert~L. Frank who attracted our attention to the papers \cite{FG} and \cite{ChSo}.
The second author is grateful to Nikolay~D. Filonov for the useful discussions. He also thanks SISSA (Trieste) 
for the hospitality during his visit in June 2013.

\end{document}